\documentclass[11pt,a4paper]{article}

\usepackage{amsmath,amssymb}
\newcommand{\e}{\varepsilon}
\newcommand{\D}{\Delta}
\newcommand{\La}{\Lambda}
\newcommand{\va}{\varphi}
\newcommand{\n}{\nabla}
\newcommand{\N}{\frac{N}{2}}

\newcommand{\g}{\int_{\mathbb{R}^{N}}}
\newcommand{\p}{\partial}

\newcommand{\R}{\mathbb{R}}
\newcommand{\h}{\hookrightarrow}

\newtheorem{definition}{Definition}
\newtheorem{theorem}{Theorem}
\newtheorem{notation}{Notation}
\newtheorem{proposition}{Proposition}
\newtheorem{corollary1}{Corollary}

\newtheorem{remarka}{Remark}

\newtheorem{lemme}{Lemma}

\title{Weak solution for compressible fluid models of Korteweg type }
\author{Boris Haspot \thanks{Universit\'e Paris XII - Val de Marne 61, avenue
du G\'en\'eral de Gaulle 94 010 CRETEIL Cedex T\'el\'ephone : (33-1)
45 17 16 51 T\'el\'ecopie : (33-1) 45 17 16 49 e-mail :
haspot@univ-paris12.fr}}
\date{}
\begin{document}
\maketitle
\subsubsection*{Abstract}
This work is devoted to proving existence of global weak solutions
for a general isothermal model of capillary fluids derived by J.E
Dunn and
J.Serrin (1985) \cite{3DS}, which can be used as a phase transition model.\\
We distinguish two cases when the dimension $N=2$ and $N=1$, in the
first case we need that $\frac{1}{\rho}\in L^{\infty}$, when $N=1$
we get a weak solution in finite time in the energy space.
\section{Introduction}
\subsection{Derivation of Korteweg model}
We are concerned with compressible fluids endowed with internal
capillarity. The model we consider  originates from the XIXth
century work by van der Waals and Korteweg \cite{3K} and was
actually derived in its modern form in the 1980s using the second
gradient theory, see for instance \cite{3JL,3TN}.
\\
Korteweg-type models are based on an extended version of
nonequilibrium thermodynamics, which assumes that the energy of the
fluid not only depends on standard variables but also on the
gradient of
the density.\\
The model derives from a Cahn-Hilliard like free energy (see the
pioneering work by J.E.Dunn and J.Serrin in \cite{3DS} and also in
\cite{3A,3C,3GP}), the conservation of mass reads:
\begin{equation}
\begin{cases}
\begin{aligned}
&\frac{\p}{\p t}\rho+{\rm div}(\rho u)=0,\\[2mm]
&\frac{\p}{\p t}(\rho u)+{\rm div}(\rho
u\otimes u)-\mu\D u-(\lambda+\mu)\n{\rm div}u+\n P(\rho)=\kappa\rho\n\D\rho,\\
\end{aligned}
\end{cases}
\label{3systeme}
\end{equation}
and:
$$\mu>0\;\;\mbox{and}\;\;2\mu+\lambda>0$$
where $\rho$ represents the density, $u\in\R^{N}$ the velocity, $P$
is a general pressure, $\mu$, $\lambda$ are the coefficients of
viscosity and $\kappa$ the coefficient of capillarity. The term
$\kappa\rho\n\D\rho$ corresponds to the capillarity
term and allows to describe the variation of density in the interfaces between two phases.\\
One can now rewrite $\rho\n\D\rho$ in the following form to
understand the difficulty of the non linear terms in distribution
sense:
\begin{equation}
\begin{aligned}
&\kappa\n\D\rho={\rm div} K,\\[2mm]
\mbox{with}\;\;K_{i,j}=\frac{\kappa}{2}&(\D\rho^{2}-|\n\rho|^{2})\delta_{i,j}-\kappa\p_{i}\rho\p_{j}\rho.\\
\end{aligned}
\end{equation}
One can recall the classical energy inequality for the Korteweg
system. Let $\bar{\rho}>0$ be a constant reference density, and let
$\Pi $ be defined
by:\\
\\
$$\Pi(s)=s\biggl(\int^{s}_{\bar{\rho}}\frac{P(z)}{z^{2}}dz-\frac{P(\bar{\rho})}{\bar{\rho}}\biggl)$$
\\
so that $P(s)=s\Pi^{'}(s)-\Pi(s)\, ,\,\Pi^{'}(\bar{\rho})=0$
and:\\
\\
$\p_{t}\Pi(\rho)+{\rm div}(u\Pi(\rho))+P(\rho){\rm div}(u)=0$ in
${\cal D}^{'}((0,T)\times\R^{N}).$
\\
\\
Notice that $\Pi$ is convex as far as $P$ is non decreasing (since
$P^{'}(s)=s\Pi^{''}(s)$), which is the case for $\gamma$-type
pressure laws or for Van der Waals law above the critical temperature.\\
Multiplying the equation of momentum conservation in the system
(\ref{3systeme}) by $\rho u$ and integrating by parts over $\R^{N}$,
we obtain the following
estimate:\\
\begin{equation}
\begin{aligned}
&\int_{\R^{N}}\big(\frac{1}{2}\rho
|u|^{2}+(\Pi(\rho)-\Pi(\bar{\rho}))+\frac{\kappa}{2}|\nabla\rho|^{2}\big)(t)dx
+2\int_{0}^{t}\int_{\R^{N}}\big(2\mu
D(u):D(u)\\
&\hspace{2cm}+(\lambda+\mu)|{\rm div} u|^{2}\big)dx
\leq\int_{\R^{N}}\big(\frac{|m_{0}|^{2}}{2\rho}+(\Pi(\rho_{0})-\Pi(\bar{\rho}))
+\frac{\kappa}{2}|\nabla\rho_{0}|^{2}\big)dx.\\
\label{3inegaliteenergie}
\end{aligned}
\end{equation}
We will note in the sequel:
\begin{equation}
{\cal E}(t)=\int_{\R^{N}}\big(\frac{1}{2}\rho
|u|^{2}+(\Pi(\rho)-\Pi(\bar{\rho}))+\frac{\kappa}{2}|\nabla\rho|^{2}\big)(t)dx,
\label{3defenergie}
\end{equation}
It follows that assuming that the initial total energy is finite:\\
$$\epsilon_{0}=\int_{\R^{N}}\big(\frac{|m_{0}|^{2}}{2\rho}+(\Pi(\rho_{0})-\Pi(\bar{\rho}))
+\frac{\kappa}{2}|\nabla\rho_{0}|^{2}\big)dx<+\infty\,,$$ then we
have the a priori bounds:
$$\Pi(\rho)-\Pi(\bar{\rho}),\;\;\mbox{and}\;\;\rho |u|^{2}\in L^{1}(0,\infty,L^{1}(\R^{N})),$$
$$\n\rho\in L^{\infty}(0,\infty,L^{2}(\R^{N}))^{N},\;\;\mbox{and}\;\;\n u\in
L^{2}(0,\infty,\R^{N})^{N^{2}}.$$ In the sequel, we aim at solving
the problem of global existence of weak solution for the system
(\ref{3systeme}). 
Assuming that we dispose from a smooth approximates sequel $(\rho_{n},u_{n})_{n\in\mathbb{N}}$ of system (\ref{3systeme}), one can remark easily that it
is difficult to pass to the limit in the quadratic term
$\n\rho_{n}\otimes\n\rho_{n}$ which belongs to $L^{\infty}(L^{1})$.
According to the classical theorems on weak topology,
$\n\rho_{n}\otimes\n\rho_{n}$ converges up to extraction to a
measure $\nu$, the difficulty is to prove that
$\nu=\n\rho\otimes\n\rho$ where $\rho$ is a limit of the sequence
$(\rho_{n})_{n\in\mathbb{N}}$
in appropriate space.\\
Another difficulty in compressible fluid mechanics is to
deal with the vacuum and we will see that this problem does appear
in the model of Korteweg, when estimating $\n\rho$. As a
matter of fact, the existence of global solution in time
for the model of Korteweg
is still an open problem.\\
The first  ones to have studied the problem, are R. Danchin and B.
Desjardins in \cite{3DD}. They showed that if we take initial data
close to a stable equilibrium in the energy space and assume that we
control the vacuum and the norm $L^{\infty}$ of the density $\rho$,
then we get some weak solution
globally in time. Controlling the vacuum amounts to bounding $\frac{1}{\rho}$ in $L^{\infty}$.\\
Recently D. Bresch, B. Desjardins and C-K. Lin in \cite{3BDL} got
some global weak solutions for the isotherm Korteweg model with some
specific viscosity coefficients. In effect, they assume that
$\mu(\rho)=C\rho$ with $C>0$ and $\lambda(\rho)=0$. In choosing
these specific coefficients they can get a gain of derivative on the
density $\rho$ and obtain an estimate for $\rho$ in $L^{2}(H^{2})$.
It is easy now with this type of estimate on $\rho$ to pass to the
limit in the term of capillarity. However a new difficulty appears
with the loss of information on the velocity $u$ and it becomes
difficult to pass to the limit in the term $\rho u\otimes u$. Hence
the solutions of D. Bresch, B. Desjardins and C-K. Lin require some
specific test
functions which depend on the density $\rho$. Indeed the loss of information is on the vacuum.\\
\\
Our result is in the same spirit as the one by Danchin and
Desjardins: we want to improve the energy inequality which allows a
gain of derivative on the density $\rho$. We show that we don't need
to control $\rho$
in $L^{\infty}$ norm to get global weak solution. However the control of the vacuum seems necessary.\\
In section \ref{3S1} we recall some definitions on the Orlicz space
and some energy inequality on the system in these spaces. In section
\ref{3S2} we show a theorem on the global existence of weak
solutions in dimension two under some conditions which amount to
controlling the vacuum. In the last section we investigate the case
of the dimension one, and we get a theorem of local existence of
solution in the energy space and a result of global existence with
small initial data.
\section{Classical a priori estimates and Orlicz spaces}
\label{3S1}
\subsection{Classical a priori estimates}
We first want to explain how it is possible to obtain natural a
priori bounds which correspond to energy when the density is close to a constant state.\\
We first rewrite the mass equation in using renormalized solutions,
and the momentum equation. We get the following formal identities:
\begin{equation}
\begin{cases}
\begin{aligned}
&\frac{1}{\gamma-1}\frac{\p}{\p
t}\big(\rho^{\gamma}-\bar{\rho}^{\gamma}
-\gamma\bar{\rho}^{\gamma-1}(\rho-\bar{\rho})\big)+{\rm
div}\big[u\frac{\gamma}{\gamma-1}(\rho^{\gamma}-\bar{\rho}^{\gamma-1}\rho)\big]=u\cdot\n(\rho^{\gamma})\\[2mm]
&\rho\frac{\p}{\p t}\frac{|u|^{2}}{2}+\rho
u\cdot\n\frac{|u|^{2}}{2}-\mu\D u\cdot u-\xi\n{\rm
div}u\cdot u+au\cdot\n\rho^{\gamma}=\kappa\rho u\cdot\n\D\rho\;,\\
\label{3a26}
\end{aligned}
\end{cases}
\end{equation}
where we note:
$\xi=\mu+\lambda$.\\
Therefore we find in summing the two equalities of (\ref{3a26}):
\begin{equation}
\begin{aligned}
\frac{\p}{\p
t}\big[\rho\frac{|u|^{2}}{2}+\frac{a}{\gamma-1}(\rho^{\gamma}+(\gamma-1)\bar{\rho}^{\gamma}
-\gamma&\bar{\rho}^{\gamma-1}\rho)\big]+{\rm
div}\big[u(\frac{a\gamma}{\gamma-1}(\rho^{\gamma}-\bar{\rho}^{\gamma-1}\rho)+\rho\frac{|u|^{2}}{2})\big]\\
&\hspace{1,8cm}-\mu\D u\cdot u-\xi\n{\rm div}u\cdot u=\kappa\rho\n\D\rho\cdot u.\\
\label{3a27}
\end{aligned}
\end{equation}
We may then integrate in space the equality (\ref{3a27}) and we get:
\begin{equation}
\begin{aligned}
&\int_{\R^{N}}\big(\rho\frac{|u|^{2}}{2}+\frac{a}{\gamma-1}(\rho^{\gamma}+(\gamma-1)\bar{\rho}^{\gamma}
-\gamma\bar{\rho}^{\gamma-1}\rho)+\kappa|\n\rho|^{2}\big)(t,x)dx+\mu\int_{0}^{t}ds\int_{\R^{N}}|D
u|^{2}dx\\
&\hspace{2cm}+\xi\int_{0}^{t}ds\int_{\R^{N}}|{\rm
div}u|^{2}dx\leq\int_{\R^{N}}\big(\rho_{0}\frac{|u_{0}|^{2}}{2}+\frac{a}{\gamma-1}(\rho_{0}^{\gamma}+(\gamma-1)\bar{\rho}^{\gamma}
-\gamma\bar{\rho}^{\gamma-1}\rho_{0})\\
&\hspace{11,7cm}+\kappa|\n\rho_{0}|^{2}\big)(x)\,dx.\\
\end{aligned}
\label{3a29}
\end{equation}
\begin{notation}
In the sequel we will note:
$$j_{\gamma}(\rho)=\rho^{\gamma}+(\gamma-1)\bar{\rho}^{\gamma}
-\gamma\bar{\rho}^{\gamma-1}\rho.$$
$${\cal E}^{\gamma}(t)=\int_{\R^{N}}\big(\rho\frac{|u|^{2}}{2}+\frac{a}{\gamma-1}(\rho^{\gamma}+(\gamma-1)\bar{\rho}^{\gamma}
-\gamma\bar{\rho}^{\gamma-1}\rho)+\kappa|\n\rho|^{2}\big)(t,x)dx$$
$${\cal E}^{\gamma}_{0}=\int_{\R^{N}}\big(\rho_{0}\frac{|u_{0}|^{2}}{2}+\frac{a}{\gamma-1}(\rho_{0}^{\gamma}+(\gamma-1)\bar{\rho}^{\gamma}
-\gamma\bar{\rho}^{\gamma-1}\rho_{0})dx$$
\end{notation}
We now want to estimate this quantity $j_{\gamma}(\rho)$ and in this
goal we recall some properties of Orlicz spaces.
\subsection{Orlicz spaces}
We begin by describing the Orlicz space in which we will work:
$$L^{q}_{p}(\R^{N})=\{f\in L^{1}_{loc}(\R^{N})/f 1_{\{|f|\leq\delta\}}\in L^{p}(\R^{N}),\;\;
f 1_{\{|f|\geq\delta\}}\in L^{q}(\R^{N})\}$$
where $\delta$ is fixed, $\delta>0$.\\
First of all, it is not difficult to check that $L^{q}_{p}$ does not
depend on the choice of $\delta>0$ since $\frac{x^{p}}{x^{q}}$ is
bounded from above and from below on any interval
$[\delta_{1},\delta_{2}]$ with $0<\delta_{1}\leq\delta_{2}<+\infty$.
In particular we deduce that we have:
$$f^{\e}\in L^{\frac{q}{\e}}_{\frac{p}{\e}}(\R^{N})\;\;\;\mbox{if}\;\;f\in L^{q}_{p}(\R^{N})\;\;\mbox{and}
\;\;p,\,q\geq\e .$$ Obviously we get
$\mbox{meas}\{|f|\geq\delta\}<+\infty$ if $f\in L^{q}_{p}(\R^{N})$
and thus we have the embedding:
$$L^{q}_{p}(\R^{N})\subset L^{q_{1}}_{p_{1}}(\R^{N})\;\;\;\mbox{if}\;\;1\leq q_{1}\leq q<+\infty,\;\;
1\leq p\leq p_{1}<+\infty. $$ Next, we choose $\Psi$ a convex
function on $[0,+\infty)$ which is equal (or equivalent) to $x^{p}$
for $x$ small and to $x^{q}$ for $x$ large, then we can define the
space $L^{q}_{p}(\R^{N})$ as follows:
\begin{definition}We define then the Orlicz space
$L^{q}_{p}(\R^{N})$ as follows:\\\texttt{}
$$L^{q}_{p}(\R^{N})=\{f\in L^{1}_{loc}(\R^{N})/\Psi(f)\in
L^{1}(\R^{N})\}.$$
\end{definition}
We can check that $L^{q}_{p}(\R^{N})$ is a linear vector space. Now
we endow $L^{q}_{p}(\R^{N})$ with  a norm so that
$L^{q}_{p}(\R^{N})$ is a separable Banach space:
$$\|f\|_{L^{q}_{p}(\R^{N})}=\inf\{t>0/\;\;\Psi(\frac{f}{t})\leq 1\}.$$
We recall now some useful properties of Orlicz spaces.
\begin{proposition}
The following properties hold:
\begin{enumerate}
\item Dual space: If $p>1$ and $q>1$
then $(L^{q}_{p}(\R^{N}))^{'}=L^{q^{'}}_{p^{'}}(\R^{N})$ where
$q^{'}=\frac{q}{q-1},\,p^{'}=\frac{p}{p-1}$.
\item $L^{q}_{p}=L^{p}+L^{q}$ if $1\leq q\leq p<+\infty$.
\item Composition: Let $F$ be a continuous function on $\R$ such that $F(0)=0$, $F$
is differentiable at $0$ and $F(t)|t|^{-\theta}\rightarrow\alpha\ne
0$ at $t\rightarrow +\infty$. Then if $q\geq\theta$,
$$F(f)\in L^{\frac{q}{\theta}}_{p}(\R^{N})\;\;\mbox{if}\;\;f\in
L^{q}_{p}(\R^{N}).$$
\end{enumerate}
\end{proposition}
Now we can recall a theorem on the Orlicz space concerning the
inequality of energy
\begin{proposition}
The function $j_{\gamma}(\rho)$ is in $ L^{1}(\R^{N})$ if and only
if $\rho-\bar{\rho}\in L^{\gamma}_{2}.$
\end{proposition}
{\bf Proof:}\\
\\
On the set $\{|\rho-\bar{\rho}|\leq\delta\}$, $\rho$ is bounded from
above, since $\gamma>1$ we thus deduce that $j_{\gamma}(\rho)$ is
equivalent to $|\rho-\bar{\rho}|^{2}$ on the set
$\{|\rho-\bar{\rho}|\leq\delta\}$.
\\
Next on the set  $\{|\rho-\bar{\rho}|\geq\delta\}$, we observe that
for some $\nu\in (0,1)$ and $C\in (1,+\infty)$, we have:
$$\nu|\rho-\bar{\rho}|^{\gamma}\leq j_{\gamma}(\rho)\leq C|\rho-\bar{\rho}|^{\gamma}.$$
\hfill {$\Box$}
\subsubsection*{Link with our energy estimate}
We recall the definition of the fractional derivative operator
$\La^{s}$:
\begin{definition}
We define the operator $\La^{s}$ as follows:
$$\widehat{\La^{s} f}=|\xi|^{s}\widehat{f}$$
\end{definition}
We recall now some useful results, we start with a proposition
coming from the theorem of interpolation by Riesz-Thorin.
\begin{proposition}
The Fourier transform is continuous from $L^{p}$ in $L^{q}$ with
$p\in[1,2]$, $q\in[2,+\infty]$ and:
$$\frac{1}{p}+\frac{1}{q}=1.$$
\end{proposition}
We recall the definition of homogeneous Sobolev space.
\begin{definition}
Let $s\in\R$. $f$ is in the homogeneous space $\dot{H}^{s}$ if:
$$|\xi|^{s}\hat{f}\in L^{2}(\R^{N}).$$
\end{definition}
\begin{proposition}
Let $f\in\dot{H}^{s}$ with $s>0$ and $f\in L^{p}+L^{2}$ with $1\leq
p<2$. Then $f\in L^{2}$. \label{3P1}
\end{proposition}
{\bf Proof:}\\
\\
Indeed we have as $f\in\dot{H}^{s}$:
$$\int_{\R^{N}}|\xi|^{2s}|\widehat{f}|^{2}d\xi<+\infty$$
so $\widehat{f}1_{\{|\widehat{f}|\geq1\}}\in L^{2}(\R^{N})$. And as
$f=f_{1}+f_{2}$ with $f_{1}\in L^{p}(\R^{N})$ and $f_{2}\in L^{2}$.
In using the Riesz-Thorin theorem, we know that
$\widehat{f_{1}}\in L^{q}(\R^{N})$ with $\frac{1}{p}+\frac{1}{q}=1$.\\
As $q\geq2$ we then have $\widehat{f}1_{\{|\widehat{f}|\leq1\}}\in
L^{2}(\R^{N})$. This concludes the proof.
\hfill {$\Box$}\\
\\
According to the above theorem and our energy estimate we get that
for all $T\in\R$,
$\rho-\bar{\rho}\in L^{\infty}(0,T;L^{\gamma}_{2}(\R^{N}))$.\\
\begin{remarka}
\label{3commentaire} We have then in using previous properties on
Orlicz spaces and (\ref{3a29}):
\begin{itemize}
\item if $\gamma\geq2$ then $L^{\gamma}_{2}(\R^{N})\hookrightarrow L^{2}(\R^{N})$ and
so $\rho-\bar{\rho}\in L^{\infty}(H^{1}(\R^{N}))$.
\item if $\gamma\leq2$ then following the proposition \ref{3P1} and the fact that $L^{\gamma}_{2}=L^{\gamma}
+L^{2}$ we get $\rho-\bar{\rho}\in L^{\infty}(H^{1}(\R^{N}))$.
\end{itemize}
\end{remarka}
We can now explain what we mean by weak solution of problem
(\ref{3systeme}) in dimension $N=1,2$.
\begin{definition}
\label{3defexistence} Let the couple $(\rho_{0},u_{0})$ satisfy;
\begin{enumerate}
\item $\rho_{0}\in L^{\gamma}_{2}(\R^{N})$, $\n\rho_{0}\in L^{2}(\R^{N})$ and
$\frac{1}{\rho_{0}}\in L^{\infty}(\R^{N})$.
\item $\rho_{0}|u_{0}|^{2}\in L^{1}(\R^{N})$
\item $\rho_{0}u_{0}=0$ whenever $x\in\{\rho_{0}=0\}$,
\end{enumerate}
We have the following definition:
\begin{itemize}
\item A couple $(\rho,u)$ is called a weak solution of problem
(\ref{3systeme}) on $I\times\R^{N}$ with $I$ an interval of $\R$ if:
\begin{itemize}
\item $\rho\in L^{\infty}(L^{\gamma}_{2}(\R^{N}))$, $\n\rho\in
L^{\infty}(L^{2}(\R^{N}))$, $\va\rho\in L^{2}(H^{1+\alpha}(\R^{N}))$
$\forall\alpha\in]0,1[$, and $\forall\va\in C^{\infty}_{0}(\R^{N})$.
\item $\frac{1}{\rho}\in L^{\infty}((0,+\infty)\times\R^{N})$.
\item $\n u\in L^{2}(L^{2}(\R^{N}))$, $\rho|u|^{2}\in
L^{\infty}(L^{1})$.
\item Mass equation holds in ${\cal D}^{'}(I\times\R^{N})$.\label{3vrai1}
\item Momentum equation holds in ${\cal D}^{'}(I\times\R^{N})^{N}$.\label{3vrai2}
\item $\lim_{t\rightarrow
0^{+}}\int_{\R^{N}}\rho(t)\va=\int_{\R^{N}}\rho_{0}\va$,
$\forall\va\in{\cal D}(\R^{N})$,\label{3vrai3}
\item $\lim_{t\rightarrow
0^{+}}\int_{\R^{N}}\rho u(t)\cdot\phi=\int_{\R^{N}}(\rho
u)_{0}\cdot\phi$, $\forall\phi\in{\cal
D}(\R^{N})^{N}$.\label{3vrai4}
\end{itemize}
\item The quantity ${\cal E}^{\gamma}_{0}$ is finite and inequality (\ref{3a29}) holds a.e in $I$.
\end{itemize}
\end{definition}
\section{Existence of global weak solutions for $N=2$}
\label{3S2}
\subsection{Gain of derivability in the case $N=2$}
In the following theorem we are interested by getting a gain of
derivative on the density $\rho$. This will enable us to treat in
distribution sense the quadratic term $\n\rho\otimes\n\rho$.
\begin{theorem}
Let $N=2$ and $(\rho,u)$ be a smooth approximate solution of the system
$(\ref{3systeme})$ such that $\frac{1}{\rho}\in
L^{\infty}((0,T)\times\R^{N})$. Then there exists a constant
$\eta>0$ depending only on the constant intervening in the Sobolev
embedding such that if:
$$\|\n\rho_{0}\|_{L^{2}(\R^{2})}+\|\sqrt{\rho_{0}}|u_{0}|\|_{L^{2}(\R^{2})}+\|j_{\gamma}(\rho_{0})\|_{L^{1}}\leq\eta$$
then we get for all $\va\in C^{\infty}_{0}(\R^{N})$:
$$\|\va\rho^{2}\|_{L^{2}_{T}(H^{1+\frac{s}{2}})}\leq M\;\;\;\;\;\mbox{with}\;\;0\leq
s<2,$$ where $M$ depends only on the initial conditions data, on
$T$, on $\va$, on $s$ and on $\|\frac{1}{\rho}\|_{L^{\infty}}$.
\label{3T1}
\end{theorem}
\begin{remarka}
In fact instead of supposing that $\frac{1}{\rho}$, we can just
assume that $\n\log\rho\in L^{\infty}(L^{2})$.
This will imply that $\rho$ will be a weight of Muckenhoupt.
\end{remarka}
\begin{remarka}
In the sequel the notation of space follows those by Runst, Sickel
in \cite{3RS}.
\end{remarka}
{\bf Proof of theorem \ref{3T1}:}\\
\\
Our goal is to get a gain of derivative on the density in using
energy inequalities and in taking advantage of the term of
capillarity. We need to localize the argument to control the low
frequencies. Let $\va\in C^{\infty}_{0}(\R^{N})$ and we have then:
\begin{equation}
\begin{aligned}
\p_{t}{\rm div}(\va\rho u)+&\p_{i,j}(\va\rho u_{i}u_{j})-(\lambda+2\mu)\D(\va{\rm div}u)+\D(\va P(\rho))\\[2mm]
&\hspace{2cm}=\frac{\kappa}{2}\D\,(\D(\va\rho^{2})-\va|\n\rho|^{2})-\kappa\p^{2}_{i,j}(\va\p_{i}(\rho)\p_{j}(\rho))+R_{\va}\\
\end{aligned}
\end{equation}
with:
$$
\begin{aligned}
&R_{\va}=\frac{\p}{\p t}(\rho u\cdot\n\va)+(\p_{i,j}\va)\rho
u_{i}u_{j}+2\p_{i}\va\,\p_{j}(\rho u_{i}u_{j})-(2\mu+\lambda)
\D\va{\rm div}u\\
&\hspace{0,5cm}-2(2\mu+\lambda)\n\va\cdot\n{\rm div}u+\D\va
a\rho^{\gamma}+2a\n\va\cdot\n(\rho^{\gamma})
-\frac{\kappa}{2}\D\va(\D\rho^{2}-|\n\rho|^{2})\\[2mm]
&\hspace{0,8cm}-\kappa\n\va\cdot\n(\D\rho^{2}-|\n\rho|^{2})+\kappa(\p^{2}_{i,j}\va)\p_{i}\rho\p_{j}\rho
+2\kappa\p_{i}\va\,\p_{j}(\p_{i}\rho\p_{j}\rho)-\kappa\D(\n\va\cdot\n\rho^{2})\\[2mm]
&\hspace{12cm}-\frac{\kappa}{2}\D(\rho^{2}\D\va).\\
\end{aligned}
$$
We can apply to the momentum equation the operator
$\Lambda^{-1}(\D)^{-1}{\rm div}$ in order to make appear a term
in $\Lambda\rho^{2}$ coming from the capillarity. Then we obtain:
\begin{equation}
\begin{aligned}
&\frac{\kappa}{2}\Lambda(\va\rho^{2})+\frac{\kappa}{2}\Lambda^{-1}(\va|\n\rho|^{2})+
\kappa\Lambda^{-1}R_{i}R_{j}(\va\p_{i}\rho\p_{j}\rho)=-\Lambda^{-3}\frac{\p}{\p t}{\rm div}(\va\rho\,u)\\[2mm]
&\hspace{0,3cm}+\Lambda^{-1}R_{i}R_{j}(\va\rho\,u_{i}u_{j})-(2\mu+\lambda)\Lambda^{-1}(\va{\rm
div}u)
+\Lambda^{-1}(\va P(\rho))+\Lambda^{-1}(\D)^{-1}R_{\va},\\
\end{aligned}
\label{3C1}
\end{equation}
where $R_{i}$ denotes the classical Riesz operator. We multiply now
the previous equality by $\La^{1+s}(\va\rho^{2})$  and we integrate
on space and in time:
\begin{equation}
\begin{aligned}
&\frac{\kappa}{2}\int^{T}_{0}\g|\Lambda^{1+\frac{s}{2}}(\va\rho^{2})|^{2}dxdt+\frac{\kappa}{2}\int^{T}_{0}
\g(|\va\n\rho|^{2})\La^{s}(\va\rho^{2})dxdt
\\[2mm]
&+\kappa\int^{T}_{0}\g\sum_{i,j}R_{i}R_{j}(\va\p_{i}\rho\p_{j}\rho)\La^{s}(\va\rho^{2})dxdt=
\g\Lambda^{-3}{\rm
div}(\va\rho\,u)\Lambda^{1+s}(\va\rho^{2})(T)dx\\
&\hspace{0,5cm}-\g\Lambda^{-3}{\rm
div}(\va\rho_{0}\,u_{0})\Lambda^{1+s}(\va\rho^{2}_{0})dx-\int^{T}_{0}\g\Lambda^{-3}{\rm
div}(\va\rho\,u)\La^{1+s}\frac{\p}{\p t}(\va\rho^{2})dxdt\\[2mm]
&+(2\mu+\lambda)\int^{T}_{0}\g\va{\rm
div}u\,\La^{s}(\va\rho^{2})dxdt-\int^{T}_{0}\g\sum_{i,j}R_{i}R_{j}(\va\rho\,u_{i}u_{j})\La^{s}(\va\rho^{2})dxdt\\[2mm]
&\hspace{2,9cm}+\int^{T}_{0}\g \va P(\rho)\La^{s}(\va\rho^{2})dxdt
+\int^{T}_{0}\g(\D)^{-1}R_{\va}\La^{s}(\va\rho^{2})dxdt.\\
\label{33}
\end{aligned}
\end{equation}
Now we want to control the term
$\int^{T}_{0}\g|\Lambda^{1+\frac{s}{2}}(\va\rho^{2})|^{2}$. Before
getting in the heart of the proof we want to rewrite the inequality
(\ref{33}) in particular the term:
$$\int^{T}_{0}\g\Lambda^{-3}{\rm
div}(\va\rho\,u)\La^{1+s}\frac{\p}{\p t}(\va\rho^{2}).$$ In this
goal we recall the renormalized equation for $\va\rho^{2}$:
\begin{equation}
\frac{\p}{\p t}(\va\rho^{2})+{\rm
div}(\va\rho^{2}u)=-\va\rho^{2}{\rm div}u+r_{\va}, \label{3C2}
\end{equation}
with $r_{\va}=-\n\va\cdot\rho^{2}u$.\\
So in using the renormalized equation (\ref{3C2}) we have:
\begin{equation}
\begin{aligned}
\int^{T}_{0}\g&\Lambda^{-3}{\rm div}(\va\rho\,u)\La^{1+s}\frac
{\p}{\p t}(\va\rho^{2})dxdt=-\int^{T}_{0}\g\Lambda^{-2}{\rm
div}(\va\rho\,u)\Lambda^{s}(\va\rho^{2}{\rm div}u)dxdt\\
&\hspace{1cm}-\int^{T}_{0}\g\Lambda^{-2}{\rm
div}(\va\rho\,u)\La^{s}{\rm div}(\va\rho^{2}u)dxdt+
\int^{T}_{0}\g\Lambda^{-3}{\rm div}(\va\rho\,u)r_{\va}.\\
\end{aligned}
\label{3C3}
\end{equation}
In combining (\ref{33}) and (\ref{3C3}) we get:
\begin{equation}
\begin{aligned}
&\frac{\kappa}{2}\int^{T}_{0}\g|\Lambda^{1+\frac{s}{2}}(\va\rho^{2})|^{2}dxdt+\frac{\kappa}{2}\int^{T}_{0}
\g(|\va\n\rho|^{2})\La^{s}(\va\rho^{2})dxdt
\\[2mm]
&+\kappa\int^{T}_{0}\g\sum_{i,j}R_{i}R_{j}(\va\p_{i}\rho\p_{j}\rho)\La^{s}(\va\rho^{2})dxdt=
\g\Lambda^{-3}{\rm
div}(\va\rho\,u)\Lambda^{1+s}(\va\rho^{2})(T)dx\\
&\hspace{0,5cm}-\g\Lambda^{-3}{\rm
div}(\va\rho_{0}\,u_{0})\Lambda^{1+s}(\va\rho^{2}_{0})dx-\int^{T}_{0}\g\Lambda^{-2}{\rm
div}(\va\rho\,u)\Lambda^{s}(\va\rho^{2}{\rm div}u)dxdt\\[2mm]
&-\int^{T}_{0}\g\Lambda^{-2}{\rm div}(\va\rho\,u) \La^{s}{\rm
div}(\va\rho^{2}u)dxdt+(2\mu+\lambda)\int^{T}_{0}\g\va{\rm
div}u\,\La^{s}(\va\rho^{2})dxdt\\[2mm]
&\hspace{1,5cm}+\int^{T}_{0}\g\sum_{i,j}R_{i}R_{j}(\va\rho\,u_{i}u_{j})\La^{s}(\va\rho^{2})dxdt+\int^{T}_{0}\g
\va P(\rho)\La^{s}(\va\rho^{2})dxdt\\[2mm]
&\hspace{1,2cm}+\int^{T}_{0}\g(\D)^{-1}R_{\va}\La^{1+\frac{s}{2}}(\va\rho^{2})dxdt
+\int^{T}_{0}\g\La^{-3}{\rm div}(\va\rho u)\La^{1+s}r_{\va}dxdt.\\
\label{3a3}
\end{aligned}
\end{equation}
\\
In order to control
$\int^{T}_{0}\g|\Lambda^{1+\frac{s}{2}}(\va\rho^{2})|^{2}$, it
suffices to bound all the other terms of (\ref{3a3}). Next we will
have
a control of $\La^{1+\frac{s}{2}}(\va\rho^{2})$ and so a gain of $\frac{s}{2}$ derivative on the density.\\
We start with the the most complicated term which requires a control
of $\frac{1}{\rho}$ in $L^{\infty}$.
\subsubsection*{1)\;$\int_{0}^{T}\g(\va|\n\rho|^{2})\La^{s}(\va\rho^{2})$:}
By induction we have $\n(\va\rho^{2})\in
L^{2}_{T}(\dot{H}^{\frac{s}{2}})$.
%
%
So by Sobolev embedding we get $\n(\va\rho^{2})\in L^{2}(L^{p})$
with
$\frac{1}{p}=\frac{1}{2}-\frac{s}{4}$ (we remark that the case $s=2$ is critical for Sobolev embedding) .\\
Now we have
$\va\n\rho=\frac{\n(\va\rho^{2})}{2\rho}-\frac{1}{2}\rho\n\va$ and
we recall that by hypothesis $\frac{1}{\rho}\in L^{\infty}$, so we
have $\va\n\rho\in L^{2}(L^{p})$ because $\rho\n\va\in
L^{\infty}(L^{r})$
for all $1\leq r\leq+\infty$ as $\rho-\bar{\rho}\in L^{\infty}(H^{1})$.\\
We now consider $\La^{s}(\va\rho^{2})$. We have by induction
$\La^{s}(\va\rho^{2})\in L^{2}(\dot{H}^{1-\frac{s}{2}})$ and
$\La^{s}(\va\rho^{2})\in L^{2}(L^{2})$ because $\va\rho^{2}\in
L^{2}(L^{2})$ which enables us to control the low frequencies of
$\La^{s}(\va\rho^{2})$. We have then $\La^{s}(\va\rho^{2})\in
L^{2}(H^{1-\frac{s}{2}})$.\\
Finally by H\"older inequality we get
$\va|\n\rho|^{2}\La^{s}(\va\rho^{2})\in L^{1}_{T}(L^{1}(\R^{N}))$
because
$\frac{1}{2}+\frac{1}{p}+\frac{1}{q}=\frac{1}{2}+\frac{1}{2}-\frac{s}{4}+\frac{s}{4}=1$
and we have:
\begin{equation}
\begin{aligned}
\int_{0}^{T}\g(\va|\n\rho|^{2})\La^{s}(\va\rho^{2})dxdt&\lesssim\|\n\rho\|_{L^{\infty}_{T}(L^{2})}
\|\La^{s}(\va\rho^{2})\|_{L^{2}_{T}(L^{q})}
\|\va\n\rho\||_{L^{2}_{T}(L^{p})},\\[2mm]
&\lesssim\|\frac{1}{\rho}\|_{L^{\infty}_{T}(L^{\infty})}
\|\Lambda^{1+\frac{s}{2}}(\va\rho^{2})\|_{L^{2}_{T}(L^{2})}^{2}\|\n\rho\|_{L^{\infty}_{T}(L^{2})}.
\end{aligned}
\end{equation}
We proceed similarly for the term:
$$\int_{0}^{T}\g\sum_{i,j}R_{i}R_{j}(\va\p_{i}\rho\p_{j}\rho)\La^{s}(\va\rho^{2})dxdt$$
because we have in following the same lines
$\va\p_{i}\rho\p_{j}\rho\in L^{2}(L^{q})$ with
$\frac{1}{q}=1-\frac{s}{4}$ and we have the Riesz operator which is
continuous from $L^{p}$ in $L^{p}$ for $1<p<+\infty$.
\\
\\
We next study the term $\g\La^{-3}{\rm
div}(\va\rho\,u)\La^{1+s}(\va\rho^{2})(t)dxdt$.
\subsubsection*{2) $\g\La^{-3}{\rm div}(\va\rho\,u)\La^{1+s}(\va\rho^{2})dx$:}
We rewrite the term $\g\La^{-3}{\rm
div}(\rho\,u)\La^{1+s}(\va\rho^{2})$ on the form:
$$
\begin{aligned}
\g\La^{-3}{\rm
div}(\va\rho\,u)\La^{1+s}(\va\rho^{2})dx=\g&\La^{-1}{\rm
div}(\va\rho\,u)\La^{-1+s}(\va\rho^{2})dx\\
&=\sum_{1\leq i\leq N}\g R_{i}(\va\rho
u_{i})\La^{-1+s}(\va\rho^{2})dx.\\
\end{aligned}
$$
As $\frac{1}{\rho}\in L_{T}^{\infty}(L^{\infty})$ then we have $u\in
L^{\infty}_{T}(L^{2})$. We recall that $\rho-\bar{\rho}\in
L^{\infty}(H^{1})$ then $\va\rho\in L^{\infty}(L^{p})$ for all
$1\leq p<+\infty$. We deduce that $\va\rho u$ belongs to
$L^{\infty}(L^{2-\beta}\cap L^{1})$ for $\beta>0$.
So we have $R_{i}(\va\rho u_{i})\in L^{\infty}_{T}(L^{r})$ for all $1<r<2$ by continuity of the operator $R_{i}$ from $L^{p}$ to
$L^{p}$ when $1<p<+\infty$.
\subsubsection*{Case $1\leq s<2$:}
Next we have:
$$\n(\va\rho^{2})=2\va\rho\n\rho+\rho^{2}\n\va$$
then we get $\n(\va\rho^{2})\in L^{\infty}(L^{2-\beta})$, in using
the fact that $\rho-\bar{\rho}\in L^{\infty}(H^{1})$ and Sobolev
embedding with H\"older inequalities.
We get then that $\va\rho^{2}$ belongs to $L^{\infty}(W^{1}_{2-\beta})=L^{\infty}(H^{1}_{2-\beta})$.\\
We have then $\La^{s-1}(\va\rho^{2})\in
L^{\infty}(H^{2-s}_{2-\beta})$. By Sobolev embedding we get
$\La^{s-1}(\va\rho^{2})\in L^{\infty}(L^{p})$ with
$\frac{1}{p}=\frac{1}{2-\beta}-\frac{2-s}{2}=\frac{1}{2-\beta}+\frac{s}{2}-1$
with $\beta$ small enough to avoid critical
embedding.\\
Finally we get $R_{i}(\va\rho u)\La^{-1+s}(\va\rho^{2})\in
L^{1}_{T}(L^{1}(\R^{N})$. Indeed we have
$\frac{1}{p}+\frac{1}{2-\beta}=\frac{2}{2-\beta}+\frac{s}{2}-1<1$ in
choosing $\beta$ small enough and
$\frac{1}{p}+\frac{1}{1+\beta}=\frac{1}{1+\beta}-1+\frac{2}{2-\beta}+\frac{s}{2}>1$
in choosing $\beta$ small enough if necessary, we conclude by
interpolation.\\
We have finally:
$$\big|\g\La^{-3}{\rm div}(\va\rho\,u)\La^{1+s}(\va\rho^{2})dx\big|\leq
M_{0}$$ with $M_{0}$ depending only on the initial data.
\subsubsection*{Case $0<s<1$:}
In this case we conclude by interpolation with the previous case.
We now want to study the other terms coming from the renormalised
equation (\ref{3C2}).
\subsubsection*{3) $\int_{0}^{T}\g\La^{-2}{\rm div}(\va\rho\,u)\La^{s}(\va\rho^{2}{\rm div}u)dxdt\;\;\mbox{and}\;\;
\int_{0}^{T}\g\La^{-2}{\rm div}(\va\rho\,u)\La^{s}({\rm
div}(\va\rho^{2}\,u))dxdt$:}
We start with:
$$\int_{0}^{T}\g\La^{-1}{\rm div}(\va\rho\,u)\La^{s-1}({\rm div}(\va\rho^{2}\,u))=
\int_{0}^{T}\g{\rm div}(\va\rho\,u)\La^{s-2}({\rm
div}(\va\rho^{2}\,u)).$$
\subsubsection*{Case $1\leq s<2$:}
We have:
$${\rm div}(\va\rho\,u)=u\cdot\n(\va\rho)+\va\rho{\rm div}u.$$
By H\"older inequalities and Sobolev embedding we get that ${\rm
div}(\va\rho\,u)$ belongs to $L^{2}_{T}(L^{2-\beta})$ for all
$\beta\in]0,1]$.
\\
Next we rewrite ${\rm div}(\va\rho^{2}\,u)$ on the form:
$${\rm div}(\va\rho^{2}\,u)=u\cdot\n(\va\rho^{2})+\va\rho^{2}{\rm
div}u.$$
As previously ${\rm div}(\va\rho^{2}\,u)$ is in
$L^{2}_{T}(L^{2-\beta})$ for all $\beta\in]0,1]$. Now by Sobolev
embedding we have $\La^{s-2}{\rm div}(\va\rho^{2}\,u)\in
L^{2}_{T}(L^{p})$ with $\frac{1}{p}=\frac{1}{2-\beta}-\frac{2-s}{2}$
with $\beta$ small enough to avoid critical Sobolev embedding.\\
We conclude that ${\rm div}(\va\rho\,u)\La^{s-2}({\rm
div}(\va\rho^{2}\,u))$ is in $L^{1}_{T}(L^{1})$ because
$\frac{1}{2-\beta}+\frac{1}{p}=\frac{2}{2-\beta}-\frac{2-s}{2}=\frac{2}{2-\beta}-1+\frac{s}{2}<1$
with $\beta$ small enough if necessary and $1+\frac{1}{p}>1$, so we
obtain the result by interpolation. Finally we have:
$$\big|\int_{0}^{T}\g\La^{-2}{\rm div}(\va\rho\,u)\La^{s}({\rm
div}(\va\rho^{2}\,u))dxdt\big|\leq M_{0}$$ with $M_{0}$ depending
only on the initial data.
\subsubsection*{Case $0<s<1$:}
We have the result by interpolation with the previous case.\\
\\
Next we proceed similarly for:
$$\int_{0}^{T}\g\La^{-1}{\rm div}(\va\rho\,u)\La^{s-1}(\va\rho^{2}{\rm div}u)dxdt.$$
\subsubsection*{4) Last terms}
We now want to concentrate us on the following term:
$$\int^{T}_{0}\g\sum_{i,j}R_{i}R_{j}(\va\rho\,u_{i}u_{j})\La^{s}(\va\rho^{2})dxdt$$
We know that $u\in L^{\infty}(L^{2})$ and $D u\in L^{2}(L^{2})$ then
$u\in L^{2}_{T}(H^{1})$ and by H\"older inequalities and Sobolev
embedding we can show that $\va\rho\,u_{i}u_{j}\in L^{2}_{T}(L^{2})$
and so $R_{i}R_{j}(\va\rho\,u_{i}u_{j})\in L^{2}_{T}(L^{2})$.\\
We have seen that $\La^{s}(\va\rho^{2})\in L^{2}(H^{1-\frac{s}{2}})$
then we have as $1-\frac{s}{2}>0$:
$$\|\La^{s}(\va\rho^{2})\|_{L^{2}_{T}(L^{2})}\leq M_{0}+\|\va\rho^{2}\|_{L^{2}_{T}(\dot{H}^{1-\frac{s}{2}})}^{\beta}$$
with $0<\beta<1$.\\
We have then:
$$|\int^{T}_{0}\g\sum_{i,j}R_{i}R_{j}(\va\rho\,u_{i}u_{j})\La^{s}(\va\rho^{2})dxdt|\leq M_{0}
+\|\La^{1+\frac{s}{2}}(\va
\rho^{2})\|^{\beta} _{L^{2}(L^{2})}$$
with $0<\beta<1$ and $M_{0}$ depending only on the initial data.\\
We are interested by the term:
$$\int_{0}^{t}\g\va{\rm div}u\,\La^{s}(\va\rho^{2})dxdt$$
We have then ${\rm div}u\in L^{2}(L^{2})$ and we have shown that
$\La^{s}(\va\rho^{2})\in L^{\infty}(L^{2})$ so we conclude in the
same way than the previous term.
\\
We finally conclude with the term:
$$\int_{0}^{T}\g \va P(\rho)\La^{s}(\va\rho^{2})dxdt.$$
Similarly we have $\La^{s}(\va\rho^{2})\in L^{2}_{T}(L^{2})$ and
$\va P(\rho)\in L^{2}_{T}(L^{2})$ because $\va\rho$ is in
$L^{\infty}(H^{1})$, and we conclude by Sobolev embedding.
\\
To finish we have to control the term with $R_{\va}$ and $r_{\va}$
that we leave to the reader. These terms are easy because they are
more regular than the preceding terms.
We finally get in using all the previous inequalities:
$$\|\La^{1+\frac{s}{2}}(\va\rho^{2})\|^{2}_{L^{2}(L^{2})}\leq C_{0}\|\n\rho\|_{L^{\infty}(L^{2})}\|\La^{1+
\frac{s}{2}}(\va\rho^{2})\|^{2}_{L^{2}(L^{2})}
+C_{1}\|\La^{1+\frac{s}{2}}(\va\rho^{2})\|^{2\beta}_{L^{2}(L^{2})}+M_{0}$$
with $0<\beta<1$ and $ C_{0}$, $C_{1}$, $M_{0}$ depends only of the norm of initial data.\\
As we have by energy inequalities
$\|\n\rho\|_{L^{\infty}(L^{2})}\leq\e<1$, we can conclude that:
$$\|\va\rho^{2}\|_{L^{2}(\dot{H}^{1+\frac{s}{2}})}\leq M_{0}$$
with $M_{0}$ depending only on the initial data. \hfill {$\Box$}
\begin{corollary1}
Let the assumptions of theorem \ref{3T1} be satisfied with the
exception of hypothesis on the viscosity coefficients which is
replaced by:
\begin{itemize}
\item it exists $c>0$, $s_{0}>0$ such that $\forall s\geq s_{0}$
$\mu(s)>c$.
\item it exists $c_{1}>0$, $m>1$ such that $\forall s\geq 0$
$\mu(s)\leq\frac{s^{m}}{c_{1}}$.
\item it exists $c_{2}>0$, $m^{'}>1$ such that $\forall s\geq 0$
$\lambda(s)\leq\frac{s^{m^{'}}}{c_{2}}$.
\end{itemize}
We get similar results as in theorem \ref{3T1} for all $\va\in C^{\infty}_{0}(\R^{N})$:
$$\|\va\rho^{2}\|_{L^{2}_{T}(\dot{H}^{1+\frac{s}{2}})}\leq M\;\;\;\;\;\mbox{with}\;\;0\leq
s<2,$$ where $M$ depends only on the initial conditions data, on
$T$, on $\va$, on $s$ and on $\|\frac{1}{\rho}\|_{L^{\infty}}$.
\label{3corollaireimport}
\end{corollary1}
{\bf Proof:}\\
\\
Indeed the only term which changes are the viscosity term.
We have then in applying the same operation as in the proof of theorem \ref{3T1} the following terms to control:
$$\int^{T}_{0}\int_{\R^{N}}\va\mu(\rho)D(u)\Lambda^{s}(\va\rho^{2})dxdt\;\;\;\mbox{and}\;\;\;\int^{T}_{0}\int_{\R^{N}}
\va R_{i,j}(\lambda(\rho)\p_{i}u_{j})\Lambda^{s}(\va\rho^{2})dxdt.$$
We conclude easily by the fact that $\sqrt{\mu(\rho)}\n u$ and $\sqrt{\lambda(\rho)}\n u$ belong to $L^{2}(L^{2})$
and as $\rho\in L^{\infty}(H^{1})$ we get that $\sqrt{\mu(\rho)}$ and $\sqrt{\lambda(\rho)}$ belongs to $L^{\infty}(L^{p})$ for all $p>2$ by Sobolev embedding.
We can then conclude by H\"older inequality.
\hfill {$\Box$}
\subsection{Existence of global weak solutions for $N=2$ away from vacuum}
We may now turn to our compactness result. First, we assume that a
sequence $(\rho_{n},u_{n})_{n\in\mathbb{N}}$ of approximate weak
solutions has been constructed by a mollifying process, which have
suitable regularity to justify the formal
estimates like the energy estimate and the previous theorem.\\
Moreover this sequence  $(\rho_{n},u_{n})_{n\in\mathbb{N}}$ has
initial data $((\rho_{0})_{n},(u_{0})_{n}))$ close to the energy
space. In using the above energy inequalities, we assume that
$j_{\gamma}((\rho_{0})_{n})$, $|\n(\rho_{0})_{n}|$ and
$(\rho_{0})_{n}|(u_{0})_{n}|^{2}$ are bounded in
$L^{1}(\R^{N})$ so that $(\rho_{0})_{n}$ is bounded in $L^{\gamma}_{2}(\R^{N})$ .\\
\\
Then it follows from the energy inequality that:
\begin{enumerate}
\item $j_{\gamma}(\rho_{n}),|\n\rho_{n}|^{2},\,\rho_{n}|u_{n}|^{2}$ are bounded in $L^{\infty}(0,T,L^{1}(\R^{N}))$,
\item $Du_{n}$ is bounded in $L^{2}(\R^{N}\times
(0,T))$,
\item $u_{n}$ is bounded
in $L^{2}(0,T,H^{1}(B_{R}))$ for all $R,T\in (0,+\infty)$.
\end{enumerate}
Extracting subsequences if necessary, we may assume that
$\rho_{n},\,u_{n}$ converge weakly respectively in
$L^{\gamma}((0,T)\times B_{R})$, $L^{2}(0,T;H^{1}(B_{R}))$ to
$\rho,\,u$ for all $R,T\in(0,+\infty)$.We also extract subsequences
for which
$\sqrt{\rho_{n}}u_{n},\,\rho_{n}u_{n},\,\rho_{n}u_{n}\otimes u_{n}$
converge weakly as previously.\\
Our goal now is to verify that the non-linear terms converge in the
sense of the distribution. The unique
difficult term to treat is $\n\rho_{n}\otimes\n\rho_{n}$.\\
Moreover we assume that $\frac{1}{\rho_{n}}$ is bounded in
$L^{\infty}(L^{\infty})$ and in using the previous theorem we get:
$$\forall\va\in C_{0}^{\infty}(\R^{N})\;\;\;\;\va\rho_{n}^{2}\;\;
\mbox{is bounded in}\;\; L^{2}(H^{1+\frac{s}{2}})\;\;\;\mbox{for all}\;\;0<s<2.$$
We can now show the following theorem:
\begin{theorem}
\label{3T2} Let $N=2$. We assume that there exists $\beta>0$ such
that for all $n\in\mathbb{N}$:
$$\rho_{n}(t,x)\geq\beta\;\;\;\mbox{for a.a}\;\;(t,x)\in (0,+\infty)\times\R^{2}$$
Then there exists $\eta>0$ such that if:
$$\|\n\rho_{0}^{n}\|_{L^{2}}+\|\sqrt{\rho_{0}^{n}}|u_{0}^{n}|\|_{L^{2}}
+\|j_{\gamma}(\rho_{0}^{n})\|_{L^{1}}\leq\eta$$ then, up to a
subsequence $(\rho_{n},u_{n})$ converges strongly to a weak solution
 $(\rho,u)$ (see definition \ref{3defexistence}) of the system (\ref{3systeme}). Moreover we
have $\n\rho_{n}\otimes\n\rho_{n}$ converges strongly in
$L^{1}_{loc}(\R\times\R^{N})$.
\end{theorem}
{\bf Proof of the theorem \ref{3T2} :}\\
\\
According to theorem \ref{3T1} we have seen that for all $\va\in
C^{\infty}_{0}(\R^{N})$ $\va(\rho_{n}^{2}-\bar{\rho}^{2})\in
L^{2}(H^{1+\frac{s}{2}})$.
\\
We can now use some results of compactness to show that
$\n\rho_{n}^{2}$ converge strongly in $L^{2}(L^{2}_{loc})$ to
$\n\rho^{2}$. We recall the following theorem from Aubin-Lions ( see
Simon for general results \cite{3Sim}).
\begin{lemme}
\label{3Aubin} Let $X\hookrightarrow Y\hookrightarrow Z$ be Hilbert
spaces such that the embedding from $X$ in $Y$ is compact. Let
$(f_{n})_{n\in\mathbb{N}}$ a sequence bounded in $L^{q}(0,T;Y)$,
(with $1<q<+\infty$) and $(\frac{d f_{n}}{dt})_{n\in\mathbb{N}}$
bounded in $L^{p}(0,T;Z)$ (with $1<p<+\infty$), then
$(f_{n})_{n\in\mathbb{N}}$ is relatively compact in $L^{q}(0,T;Y)$.
\end{lemme}
We need to localize because we have some result of compactness for
the local Sobolev space $H_{loc}^{\frac{s}{2}}$ which is compactly
embedded in $L^{2}_{loc}$. Let $(\chi_{p})_{p\in\mathbb{N}}$ be a
sequence of $C^{\infty}_{0}(\R^{N})$ cut-off functions supported in
the ball $B(0,p+1)$ of $\R^{N}$ and equal to 1 in a neighborhood of
$B(0,p)$.\\
We have then in using mass equation:
$$\frac{d}{dt}\n(\rho_{n}^{2})+\n{\rm
div}(\rho^{2}_{n}u_{n})=-\n(\rho_{n}^{2}{\rm div}u_{n})$$ We can
then show that
$\big(\frac{d}{dt}(\chi_{p}\n(\rho_{n}^{2}))\big)_{n\in\mathbb{N}}$
is uniformly bounded for all $p$ in $L^{q}(H^{\alpha})$ for
$\alpha<0$ in using energy inequalities and
$\big(\chi_{p}\n(\rho_{n}^{2})\big)_{n\in\mathbb{N}}$ is uniformly
bounded for all $p$ in $L^{2}(H^{\frac{s}{2}})$. Apply lemma
\ref{3Aubin} with the family
$\big(\chi_{p}\n(\rho_{n}^{2})\big)_{n\in\mathbb{N}}$ and
$X=\chi_{p}H^{\frac{s}{2}}$, $Y=\chi_{p}L^{2}$, $Z=\chi_{p}
H^{\alpha}$, then use Cantor's diagonal process. This finally
provides that:
\begin{equation}
\forall
p>0\;\;\;\chi_{p}\n(\rho_{n}^{2})\rightarrow_{n\rightarrow+\infty}
\chi_{p}\n(\rho^{2})\;\; \mbox{in}\;\;L^{2}(L^{2}). \label{3egalite}
\end{equation}
We now want to show some results of strong convergence for $\n\rho_{n}$.\\
We have then for all $\phi\in C_{0}^{\infty}$:
$$|\phi(\n\rho_{n}-\n\rho)|\leq\frac{1}{\rho_{n}}|\phi(\n\rho_{n}^{2}-\n\rho^{2})|+|\n\rho^{2}|
|\phi(\frac{1}{\rho_{n}}-\frac{1}{\rho})|=A_{n}+B_{n}.$$ We have
$A_{n}$ converges to $0$ in $L^{2}(L^{2})$ because
$\frac{1}{\rho_{n}}$ is uniformly bounded in
$L^{\infty}(L^{\infty})$ and in using (\ref{3egalite}) for all
$\phi\in C_{0}^{\infty}$, $\phi\n\rho_{n}^{2}$ converges to
$\n\rho_{n}^{2}$ in $L^{2}(L^{2})$. In the same way we have
$\frac{1}{\rho_{n}}$ converges a.e to $\frac{1}{\rho}$ (because we
can show in using again the theorem of Aubin-Lions that
$\rho_{n}-\bar{\rho}\rightarrow\rho-\bar{\rho}$ in $L^{2}(L^{q})$ up
to an extraction with $1\leq q<+\infty$ then we can extract again a
subsequence so that $\rho_{n}$ converges a.e to $\rho$). Moreover we
have:
$$\big|\frac{1}{\rho_{n}}-\frac{1}{\rho}\big|\leq\frac{2}{\beta}$$
then by the theorem of dominated convergence we have $B_{n}$ tends to $0$ in  $L^{2}(L^{2})$.\\
We can conclude that:
$$\forall\phi\in C_{0}^{\infty},\;\;\;\phi\n\rho_{n}\otimes\n\rho_{n}\rightarrow_{n}\phi\n\rho\otimes\n\rho
\;\;\;\mbox{in}\;\;L^{1}(L^{1}).$$ \hfill {$\Box$}
\section{Existence of weak solution in the case $N=1$}
\label{3S3} We are now interested by the case $N=1$. To start with,
we focus on the gain of derivative for $\n\rho^{2}$.
\subsection{Gain of derivative}
We can now write a theorem where we expose a gain of derivative on
the density $\rho$ in using the same type of inequalities as in the
case $N=2$.
\begin{theorem}
Let $(\rho,u)$ be a regular solution of the system (\ref{3systeme})
with initial data in the energy space. Then we have:
$$\|\rho^{2}\|_{L^{2}_{T}(H^{1+\frac{s}{2}}(\R))}\leq M_{0}$$
with $0\leq s<\frac{1}{2}$ and $M_{0}$ depending only of the initial
data. \label{3T3}
\end{theorem}
\begin{remarka}
We observe the two important facts:
\begin{enumerate}
\item We don't need any hypothesis on the size of the initial
data.
\item We don't need to localize because we know that $\rho\in L^{\infty}_{t,x}$.
\item We don't need to assume that $\frac{1}{\rho}\in L^{\infty}$, this is important.
In fact if we assume that $\rho_{0}$
admits some vacuum, we could show that
$\n\rho_{n}(t,x)\rightarrow\n\rho$ a.e on the set $A=\{\rho>0\}$
because we have by compactness up to an extraction
$\n\rho_{n}^{2}\rightarrow\n\rho^{2}$ a.e.
\end{enumerate}
\end{remarka}
{\bf Proof of theorem \ref{3T3} :}\\
\\
We use the same estimates as in the previous proof except for the
delicate term:
$\int_{0}^{T}\int_{\R}|\p_{x}\rho|^{2}\La^{s}\rho^{2}$.\\
We have then $\p_{x}\rho\in L^{\infty}(L^{2})$ and
$\rho-\bar{\rho}\in L^{\infty}(L^{2})$ so $\rho-\bar{\rho}\in
L^{\infty}(H^{1})$ and we have then by Sobolev embedding $\rho\in
L^{\infty}(L^{\infty})$. By composition theorem on the Sobolev space
we get $(\rho-\bar{\rho})^{2}\in L^{\infty}(H^{1})$ and
$\rho^{2}-\bar{\rho}^{2}=(\rho-\bar{\rho})^{2}+2\bar{\rho}(\rho-\bar{\rho})\in
L^{\infty}(H^{1})$. Finally we get for $0<s\leq1$,
$\La^{s}\rho^{2}\in L^{\infty}(H^{1-s})$. Now for $0\leq
s<\frac{1}{2}$ by Sobolev embedding we obtain:
$$\La^{s}\rho^{2}\in L^{\infty}(L^{\infty}).$$
So we can control the term
$\int_{0}^{T}\int_{\R}|\p_{x}\rho|^{2}|\La^{s}\rho^{2}$ as follows:
$$\int_{0}^{T}\int_{\R}|\p_{x}\rho|^{2}|\La^{s}\rho^{2}|\lesssim\|\p_{x}\rho\|_{L^{\infty}_{T}(L^{2})}^{4}.$$
We treat the other terms similarly as in the previous proof. \hfill{$\Box$}
\subsection{Results of compactness}
We can now prove our result of stability of solution in the case
$N=1$ in using the previous gain of derivative. Let
$(\rho_{n},u_{n})_{n\in\mathbb{N}}$ a sequel of approximate weak
solutions of system (\ref{3systeme}).
\begin{theorem}
Let $(\rho_{0}^{n},u_{0}^{n})$ initial data of the system
(\ref{3systeme}) in the energy space what it means that:
$$\int_{\R}\rho_{0}^{n}|u_{0}^{n}|^{2}+j_{\gamma}(\rho_{0}^{n})+|\p_{x}\rho_{0}^{n}|^{2}dx\leq M$$
with $M>0$.\\
Moreover we assume that $\rho_{0}^{n}\geq c>0$. Then there exists a
time $T$ such that up to a subsequence, $(\rho_{n},u_{n})$ converges
strongly to a weak solution $(\rho,u)$
on $(0,T)\times\R$ in the sense of the distribution (see definition \ref{3defexistence}).\\
Moreover $\p_{x}\rho_{n}$ converges strongly in
$L^{2}(\R\times\R)$ to $\p_{x}\rho$. \label{3T4}
\end{theorem}
{\bf Proof of the theorem \ref{3T4}:}\\
\\
We want now to control the vacuum of $\frac{1}{\rho}$. We recall
that $\rho-\bar{\rho}\in L^{\infty}(H^{1})$ and $\frac{\p}{\p t}\rho\in
L^{\infty}(L^{p})$, then we have by Aubin-Lions theorem for all
$\phi\in C_{0}^{\infty}$, $\phi(\rho-\bar{\rho})\in
C([0,T],L^{\infty})$.
\\
Moreover we set:
$$f_{n}(t)=\int_{\R}|\p_{x}\rho|^{2}1_{\mathbb{R}\setminus B(0,n)}$$ then we have $f_{n}$ is
a decreasing sequence and converges to $0$. Then in using the
theorem of Dini we get that on $[0,T]$,  $f_{n}$ converges uniformly
to $0$. Let $\e>0$, then it exists $n_{0}\in\mathbb{N}$ enough big such that for all $n\geq n_{0}$:
$$\int_{\R}|\p_{x}(\rho-\bar{\rho})|^{2}1_{\mathbb{R}\setminus B(0,n)}\leq\e$$
we proceed similarly to show that for a $n^{'}_{0}$ big enough we
have for all $n^{'}\geq n^{'}_{0}$:
$$\int_{\R}|\rho-\bar{\rho}|^{2}1_{\mathbb{R}\setminus B(0,n^{'})}\leq\e$$
then we have for $\phi\in C_{0}^{\infty}$ with a big enough support
$(\rho-\bar{\rho})(1-\phi)\in L^{\infty}(H^{1})$ with a norm
inferior to $\e$ then it exists $c$ such that by continuity with
$0<\e\leq c<\bar{\rho}$, $|\rho-\bar{\rho}|\leq c$.
\\
Finally we have $\frac{1}{\rho}\in L^{\infty}(L^{\infty})$ and we
can conclude as in the previous theorem \ref{3T3} show that:
$$\forall \phi\in C^{\infty}_{0}\;\;\;\phi\p_{x}\rho_{n}\rightarrow_{n}\phi\p_{x}\rho\;\;\;\mbox{in}\;\;L^{2}(L^{2}).$$
\hfill {$\Box$}
\begin{theorem}
Let $(\rho_{0}^{n},u_{0}^{n})$ initial data of the system
(\ref{3systeme}) in the energy space.\\
Then it exists $\e>0$ such that if:
$$\|\p_{x}\rho_{0}^{n}\|_{L^{2}}+\|\sqrt{\rho_{0}^{n}}|u_{0}^{n}|\|_{L^{2}}
+\|j_{\gamma}(\rho_{0}^{n})\|_{L^{1}}\leq\e$$
then up to a subsequence $(\rho_{n},u_{n})$ converges strongly to a
weak solution  $(\rho,u)$ on $\R\times\R$ (see the definition
\ref{3defexistence}).\\
Moreover $\p_{x}\rho_{n}$ converges strongly to $\p_{x}\rho$ in
$L^{2}(\R\times\R)$. \label{3T5}
\end{theorem}
{\bf Proof of theorem \ref{3T5}:}\\
\\
We can show easily that
$\|\rho_{n}-\bar{\rho}\|_{L^{\infty}(H^{1})}\leq C\e$ and so we have
$\frac{1}{\rho_{n}}$ is uniformly bounded in
$L^{\infty}(L^{\infty})$. We can then conclude as in the proof of
theorem
\ref{3T3}.\\
\null\hfill {$\Box$}

\end{document}